\documentclass[oneside]{amsart}

\usepackage{graphicx}

\usepackage{textcomp}

\newtheorem{theorem}{Theorem}[section]
\newtheorem{lemma}[theorem]{Lemma}
\newtheorem{proposition}[theorem]{Proposition}
\newtheorem{corollary}[theorem]{Corollary}

\theoremstyle{definition}
\newtheorem{definition}[theorem]{Definition}

\theoremstyle{remark}
\newtheorem{remark}[theorem]{Remark}

\numberwithin{equation}{section}

\newcommand{\majN}{\mathcal{N}}
\newcommand{\majL}{\mathcal{L}}
\newcommand{\ind}[1]{\mathbf{1}_{#1}}
\newcommand{\limNT}[2]{{\displaystyle \lim_{\stackrel{N.T.}{#1\to#2}}}}
\newcommand{\simsubset}{\stackrel{\sim}{\subset}}

\begin{document}

\title{Harmonic functions on hyperbolic graphs}

\author{Camille PETIT}
\address{Universit\'e Joseph Fourier Grenoble 1\\
Institut Fourier UMR 5582 UJF-CNRS\\
100 rue des Maths, BP 74\\
38402 Saint Martin d'H\`eres\\
France }
\email{camille.petit@ujf-grenoble.fr}


\date{2010}



\begin{abstract}
We consider admissible random walks on hyperbolic graphs. For a given harmonic function on such a graph, we prove that asymptotic properties of non-tangential boundedness and non-tangential convergence are almost everywhere equivalent. The proof is inspired by the works of F. Mouton in the cases of Riemannian manifolds of pinched negative curvature and infinite trees. It involves geometric and probabilitistic methods.

\end{abstract}

\maketitle

\tableofcontents


\section{Introduction}

The study of non-tangential convergence of harmonic functions began in 1906 with P.~Fatou \cite{Fat06}, who showed that a given positive harmonic function on the unit disc of $\mathbb{R}^2$ admits at almost all points of the unit circle a non-tangential limit. The same is true in many general cases: euclidean half-spaces, trees (\cite{Car72}), free groups (\cite{Der75}), Riemannian manifolds of pinched negative curvature (\cite{AnS85}, \cite{Anc87}) and Gromov hyperbolic graphs (\cite{Anc90}). One is thus naturaly led to the study of cases where the harmonic function is not necessarily positive. Fatou's conclusion is no longer true in this more general case, and several authors have made attempts to give criteria for the harmonic function to admit non-tangential limit at a point of the boundary. In the case of the euclidean half space ${\mathbb{R}^n\times\mathbb{R}_+^*}$, A.P.~Calder\`{o}n and E.M.~Stein (\cite{Cal50a}, \cite{Cal50b}, \cite{Ste61}) proved that for a harmonic function $u$, the following three properties are equivalent for almost all point $\theta$ of the boundary:
\begin{itemize}
	\item the function $u$ is non-tangentially convergent at $\theta$
	\item the function $u$ is non-tangentially bounded at $\theta$
	\item the area integral $\intop_{\Gamma_\theta} |\nabla u(x,y) |^2y^{1-n}dxdy$ is finite (for all $\Gamma^\theta$ where $\Gamma^\theta$ is a non-tangential cone).
\end{itemize}	
In 1978, using probabilistic methods, J.~Brossard proved the same result \cite{Bro78}. Shortly after, A.~Koranyi remarked that hyperbolic spaces provide a more natural setting for this study. Indeed, several notions have simpler expressions in this case. For instance the boundary becomes an ideal one, non-tangential cones become tubular neighborhoods of geodesic rays. Following this remark, F.~Mouton proved in 1994 an analogous result for harmonic functions on Riemannian manifolds of pinched negative curvature  \cite{Mou94}, and in 2000 for harmonic functions on trees \cite{Mou00}. We prove here a partial analogue for hyperbolic graphs:
for a harmonic function (the notion of harmonicity is here relative to a random walk on the graph), non-tangential convergence is almost everywhere equivalent to non-tangential boundedness.

We introduce in section~2 the notions of random walks and harmonic functions on hyperbolic graphs and in section~3 the boundary at infinity, which enables us to state our main result in section~4. Section~5 is devoted to the conditioning of the random walk to exit at a fixed point of the boundary and to the proof of a stochastic result. In order to prove the non-tangential convergence criterion, we state geometric lemmas in section~6. We then prove the main result in section~7.\\


\section{Harmonic functions on hyperbolic graphs}

We shall briefly introduce the notions of hyperbolic graphs, random walks, harmonic functions and Green functions. The reader can refer to \cite{Woe00} for more details.

\subsection{Hyperbolic graphs}

Gromov hyperbolicity was introduced in the 80's by M.~Gromov \cite{Gro81}. One way to define it is the following:

\begin{definition}
On a metric space $(X,d)$, one defines the Gromov product of two points $x,y\in X$ with respect to $o\in X$ by
$$(x,y)_o=\frac{1}{2}[d(x,o)+d(y,o)-d(x,y)].$$
For a real $\delta \geq 0$, a metric space $X$ is said to be $\delta$-\textit{hyperbolic} if for all $x,y,z,o\in X$, 
$$(x,z)_o\geq \min \{ (x,y)_o,(y,z)_o \} -\delta.$$
\end{definition} 

A metric space $(X,d)$ is geodesic if for every pair of points $x$ and $y$ in $X$, there is a geodesic segment (not necessarily unique) joining $x$ to $y$ in $X$ ( \textit{i.e.} an isometric embedding of the real interval $[0,d(x,y)]$ into $X$ which sends $0$ to $x$ and $d(x,y)$ to $y$).

The definition of Gromov hyperbolicity makes sense in all metric spaces. However, it has a nice geometric interpretation when the space is geodesic. A geodesic triangle consists of three points $x,y,z\in X$ together with geodesic segments $\alpha$, $\beta$, $\gamma$ (respectively from $y$ to $z$, $z$ to $x$ and $x$ to $y$) called the sides. A triangle is called $\eta$-thin for a real $\eta\geq 0$ if every point of a side is at distance at most $\eta$ from the union of the other two sides. If a geodesic metric space $X$ is $\delta$-hyperbolic, then every geodesic triangle in $X$ is $4\delta$-thin. Remark that the converse also holds, if every geodesic triangle in $X$ is $\eta$-thin, then $X$ is $3\eta$-hyperbolic. 
The reader can keep in mind that the Gromov product $(x,y)_o$ can be seen as a rough measure of the distance between $o$ and a geodesic segment joining $x$ and $y$ (see \cite{GdlH90}). Precisely, if $X$ is $\delta$-hyperbolic and $\gamma$ is a geodesic segment from $x$ to $y$, then
\begin{equation}
d(o,\gamma)-2\delta \leq (x,y)_o \leq d(o,\gamma).   \label{thin triangles}
\end{equation}

Let $S$ be a \textit{countable graph}, that is, it is a countable set $S$ equipped with a reflexive symmetric relation $Z\subset S\times S$. A path from $x$ to $y$ in $S$ is a sequence $[x=x_{0},x_{1},...,x_{k}=y]$
such that for all indices $i$, $(x_{i-1}, x_{i})\in Z$. The integer $k$ is
the length of the path. We shall always assume that $S$ is \textit{connected}, \textit{i.e.} that for every pair $x$, $y$ in $S$, there is a path from $x$ to $y$. The graph $S$ carries an integer-valued metric: $d(x,y)$ is the minimum among all the lengths of the paths from $x$ to $y$.
If the metric space $(S,d)$ is $\delta$-hyperbolic for a real $\delta\geq 0$, we will say that the graph $S$ is hyperbolic.
Typical examples are provided by trees and Cayley graphs of certain groups. If the Cayley graph of a finitely generated group associated to one finite generating set is hyperbolic, then the Cayley graph associated to every finite generating set is hyperbolic and the group is called Gromov hyperbolic. 

We will focus our interest on coercive graphs satisfying the geometric assumption (GRBD).

\begin{definition}
A graph $S$ is called \textit{coercive} if there is some positive $\alpha$ such that the following Poincar\'e-Sobolev inequality holds:
$$\sum_{(x,y)\in Z} |u(y)-u(x)|^2 \geq \alpha \cdot ||u||_2^2$$
for all $u:S\to\mathbb{R}$ with finite support.
\end{definition}

One can verify that the coercivity of $S$ is equivalent to an isoperimetric inequality \cite{Anc90}. An example of coecive graph is given  by a discrete approximation of the hyperbolic ball: denote by $B$ the unit ball of $\mathbb{R}^n$ equipped with the hyperbolic metric $d_h$. Let $S\subset B$ be such that for all $x\in B$, $d_h(x,S)\leq c_1$ and for all $x,y\in S$, $d_h(x,y)\geq c_2$. Equip $S$ with the relation $(x,y)\in Z$ iff $d_h(x,y)\leq 3c_1$. Then, $S$ is a coercive graph, whose metric is uniformly equivalent to $d_h$.

\begin{definition}
A graph $S$ satisfies the \textit{Geodesic Ray at Bounded Distance (GRBD)} assumption if, given a base point $o$, there is a constant $K\geq 0$ with the property that every point in the graph is at distance at most $K$ from a geodesic ray starting from $o$.
\end{definition}

This assumption will be used in the proof of lemma \ref{lem:geometric lemma2}. It is satisfied for instance by geodesically complete graphs (any two points in the graph can be joined by a geodesic line) and by Cayley graphs of hyperbolic groups. Indeed, assume that $S$ is a Cayley graph of a hyperbolic group and denote by $\delta$ a hyperbolicity constant of $S$. Let $x\in S$. Choose two arbitrary points $\xi_1,\xi_2\in\partial S$ and a geodesic joining them, which exists by visibility (see \cite{GdlH90}). Up to translation, one can assume without loss of generality that it contains $x$.
Let $\gamma_1$ denote a geodesic ray from $o$ to $\xi_1$ and $\gamma_2$ a geodesic ray from $o$ to $\xi_2$. Because all triangles are $4\delta$-thin, $x$ is at distance at most $4\delta$ from one of the two geodesic rays.

\subsection{Random walks}

Let $S$ be a graph and $d$ the corresponding distance.
Let us choose a (transition) function $p:S\times S\to \mathbb{R}^+$. This function $p$ is \textit{markovian} (resp. \textit{submarkovian}) if for all $x\in S$, $\sum_{y\in X} p(x,y)= 1$ (resp. $\sum_{y\in X} p(x,y) \leq 1$)  and \textit{admissible} in the sense of A.~Ancona (\cite{Anc88}) if the following relations hold:
\begin{enumerate}
	\item $\exists c_0>0, \exists \ell\in\mathbb{N}^*$  such that  $\forall x,y\in S$, $d(x,y)\leq 1 \Rightarrow \sum_{0\leq j\leq \ell} p^j(x,y)\geq c_0$,
	\item $\exists m_1\in\mathbb{N}^*$ such that  $\forall x,y\in S$, $p(x,y)>0 \Rightarrow d(x,y)\leq m_1$.
\end{enumerate}	

The admissible conditions are geometric adaptedness properties of the transition function $p$ to the structure of the graph $S$.

Remark that the adjoint kernel $p^*(x,y)=p(y,x)$ and $p+tI$, $t\geq0$ are also admissible.
In the following, we will always assume $p$ to be markovian. We define the $p$-random walk on $S$ as the Markov chain with state space $S$ and transition probabilities $p(x,y)$, $x,y\in S$.
It is given by a family of random variables $(X_n)_{n\in \mathbb{N}}$ on $S$ where $X_n$ is the position at time $n$. We can choose the probability space to be the space $\Omega=\mathcal{C}(\mathbb{N},S)$ of all infinite paths (then, $X_n(\omega)=\omega(n)$), equipped with the $\sigma$-algebra arising from the countable product of $\mathcal{P}(S)$. We will denote by $(\mathbb{P}_z)_{z\in S}$ the law of this random walk, where $\mathbb{P}_z$ is the probability obtained when the walk starts from $z$ and by $\mathcal{F}_n$ the $\sigma$-algebra generated by $X_i,i\leq n$.

We can now state a classical property which will be useful in the following. For an almost surely finite stopping time $T$, we denote by $\Theta^T$ the map: $\Theta^T(\omega)=\omega( \cdot +T(\omega))$.

\begin{lemma}[Strong Markov property]   \label{lem:Strong Markov roperty}
For a non-negative random variable $F$ on $\Omega$ and an almost surely finite stopping time $T$ one has
$$\mathbb{E}_x[ F\circ \Theta^T | \mathcal{F}_T] =u_F (X_T) \text{ where } u_F(y)=\mathbb{E}_y[F]. $$
\end{lemma}

\subsection{Harmonic functions and the Green function}

We associate to the random walk a Laplace operator $\Delta$ which
acts on functions $f:S\to\mathbb{R}$ by
\[
\Delta f(x)=\mathbb{E}_{x}[f(X_{1})]-f(x)=\sum_{y\in S}p(x,y)f(y)-f(x).
\]
A function $f$ is said to be \textit{harmonic} if $\Delta f=0$ and \textit{superharmonic} if $\Delta f\leq 0$.

The \textit{Green function} associated to the random walk is thus defined on $S\times S$ by 
$$G(x,y):=\sum_{n=0}^{\infty} \mathbb{P}_{x}[X_{n}=y]=\mathbb{E}_{x}\left[ \sum_{n=0}^{\infty}\ind{\{ X_{n}=y\} }\right].$$
It can be seen by the Markov property that the function $G(\cdot ,y)$ is harmonic on $S\setminus\{y\}$ and superharmonic on $S$.

\begin{lemma}[Martingale property]  \label{lem:Martingale property}  
Let $z$ be a point of $S$ and $f$ be a function on $S$. Then, the sequence of random variables
$$M_n=f(X_n)-\sum_{k=0}^{n-1}\Delta f(X_k)$$
is a $(\mathcal{F}_n)$-martingale for the probability $\mathbb{P}_z$. In particular, $(f(X_n))_n$ is a martingale if $f$ is harmonic.
\end{lemma}

Denote by $G^t$ the Green kernel of the admissible transition function $p+tI$. A.~Ancona introduced in \cite{Anc88} the condition:
$$\text{(*) there exists } \varepsilon>0, \text{ such that } G^{\varepsilon} \text{ is finite.}$$

This hypothesis implies in particular that the random walk is transient and that some Harnack inequalities hold at infinity. In \cite{Anc90}, A.~Ancona proves the following proposition:

\begin{proposition}[Ancona]
\label{prop:exponential decay}
If $S$ is coercive, every admissible kernel $p$ on $S$ such that $p$ and $p^*$ are submarkovian satisfies condition (*) and admits a Green function $G$ such that $G(x,y)\leq C \cdot \exp(-\beta d(x,y))$ for some positive constants $C,\beta$.
\end{proposition} 

\begin{remark}
 the assumption in the proposition is satisfied by non-amenable, finitely generated groups $\Gamma$ with an admissible probability measure $\nu$ on $\Gamma$ ($p(x,y)=\nu(\{ x^{-1}y\} )$).
\end{remark}


\section{Boundary at infinity}

Since we are interested in non-tangential convergence, we need a notion of boundary. In fact, we will focus our interests on two types of boundaries: the geometric boundary and the Martin one. To define them, we will need to fix a base point $o\in S$, but the compactifications below do not depend on the choice of $o$.
\begin{itemize}
	\item The \textit{geometric boundary $\partial S$}. Assume $S$ to be a hyperbolic graph. Let us denote by $E$ the set of sequences $(x_i)_i$ in $S$ such that $\lim_{i,j\to\infty} (x_i,x_j)_o=+\infty$ and by $\partial S=E/\mathcal{R}$ the set obtained by factoring $E$ with respect to the equivalence relation: $(x_i)_i \mathcal{R} (y_j)_j$ iff $\lim_{i,j\to\infty } (x_i,y_j)_o=+\infty$.  An equivalent way to describe $\partial S$ is via equivalence of geodesic rays (see \cite{GdlH90}): two geodesic rays $\gamma_1$ and $\gamma_2$ are equivalent if $\liminf_{k\to\infty} d(\gamma_1(k),\gamma_2(\mathbb{N}))<+\infty$.
	One can extend the Gromov product to two points $x,y\in\partial S$ (resp. $x\in S$, $y\in\partial S$ or $x\in\partial S$, $y\in S$) with 
	$$(x,y)_o=\sup \liminf_{i,j\to\infty} (x_i,y_j)_o  \hspace{0.5cm} ( \text{resp. } (x,y)_o=\sup \liminf_{j\to\infty} (x,y_j)_o ),$$
	where the supremum is taken over all sequences $(x_i)_i$ in the class of $x$ and $(y_j)_j$ in the class of $y$.
	For a real $r>0$ and a point $x\in\partial S$, denote by $V_r(x)=\{ y\in S\cup \partial S \, | \, (x,y)_o \geq r \}$. We equip $S\cup\partial S$ with the unique topology containing open sets of $S$ and admitting the sets $V_r(x)$ with $r\in \mathbb{Q}^+$ as neighborhood base at any $x\in\partial S$.
	It provides a compactification $\tilde S$ of $S$ (that is a compact Hausdorff space with countable base of the topology such that $S$ is open and dense in $\tilde S$). The compactification $\tilde S$ can also be obtained as the completion of $S$ for a good choice of a metric on $S$ (see \cite{Woe00}).

	\item The \textit{Martin boundary}. Assume the random walk to be transient. One defines the Martin kernel by $K(x,y)=\frac{G(x,y)}{G(o,y)}$. The Martin compactification $\hat S$ is the unique smallest compactification of $S$ for which all kernels $K(x,.)$, $x\in S$, extend continuously. The Martin boundary is $\hat S\setminus S$. A sequence $(y_i)_i\in S^{\mathbb{N}}$ converges to the Martin boundary if $d(o,y_i)\to\infty$ and $(K(.,y_i))_i$ converges pointwise. Two such sequences are equivalent if their limits coincide at each point of $S$. The Martin boundary allows us to represent non-negative harmonic functions by non-negative measures on this boundary (see \cite{Woe00}).
	
\end{itemize}	

Results by A.~Ancona (\cite{Anc90}) imply that in the case of a hyperbolic graph with an admissible transition function $p$ on $S$ satisfying condition (*), these two compactifications coincide. In the following, we will assume that $S$ is a coercive hyperbolic graph and $p$ be an admissible markovian transition function on $S$ such that $p^*$ is submarkovian. The above two boundaries of $S$ thus coincide and we shall denote it by $\partial S$. There is a $\partial S$-valued random variable $X_\infty$ such that the random walk $(X_n)_n$ converges $\mathbb{P}_z$-almost surely to $X_\infty$ for all $z\in S$ (see \cite{Woe00}). When dealing with harmonic functions and random walks, there is a natural family of measures on $\partial S$ called the harmonic measures $\mu_z, z\in S$. The measure $\mu_z$ is the distribution of the random variable $X_\infty$ when the walk starts from $z$. Different measures $\mu_z$ are equivalent, so we can define a notion of $\mu$-negligeability. Their Radon-Nykodim derivatives are given by
$$ \frac{d\mu_y}{d\mu_x}(\theta)=\lim_{z\to\theta} \frac{G(y,z)}{G(x,z)}.$$
These measures allow to represent bounded harmonic functions by the Poisson formula (see \cite{Woe00} and lemma \ref{lem:Representation of bounded harmonic functions}).


\section{Main result}

\textbf{Setting:} We fix now a coercive hyperbolic graph $S$ satisfying (GRBD) and an admissible markovian transition function $p$ on $S$ such that $p^*$ is submarkovian.

Let $d$ be the canonical distance on $S$ , $\delta\geq 0$ the hyperbolicity constant, $c_0$, $\ell$ and $m_1$ the admissibility constants and $o\in S$ a base point. After possibly enlarging it, we will assume that $\delta$ is an integer strictly bigger than 3.

Let us now define the non-tangential notions. If $c>0$ and $\theta\in\partial S$, let us denote by
$$\Gamma_c^\theta :=\{ x\in S \, | \, \exists \gamma \text{ a geodesic ray from } o \text{ to } \theta \text{ such that } d(x,\gamma)<c \}$$
the non-tangential tube of radius $c$ and vertex $\theta$. A function $u$ \textit{converges non-tangentially} at $\theta$ if, for all $c>0$, $u(x)$ has a limit as $x$ goes to $\theta$ in $\Gamma_c^\theta$. In the same manner, the function $u$ is \textit{non-tangentially bounded} at $\theta$ if, for all $c>0$, $u$ is bounded on $\Gamma_c^\theta$.
Remark that these non-tangential notions do not depend on $o$, due to the alternative definition of $\partial S$ by geodesic rays.

We can now state our main result.

\begin{theorem}    \label{thm:main result}
In the setting above, for a harmonic function $u$, the following two properties are equivalent for $\mu$-almost all $\theta\in\partial S$:
\begin{enumerate}
	\item the function $u$ converges non-tangentially at $\theta$,
	\item the function $u$ is non-tangentially bounded at $\theta$.
\end{enumerate}
\end{theorem}	

Denote 
$$\majL_c=\{ \theta\in\partial S \, | \, \lim_{\stackrel{x\in\Gamma_c^\theta}{x\to\theta}} u(x) \text{ exists and is finite} \},$$
$$\majN_c=\{ \theta\in\partial S \, | \, N_c^\theta (u)<\infty \} \text{ where } N_c^\theta (u)=\sup_{x\in\Gamma_c^\theta} |u(x)|$$
and observe that
$$\majL=\bigcap_{c>0} \majL_c   \hspace{0,5cm} \text{ and } \hspace{0,5cm}  \majN=\bigcap_{c>0}\majN_c.$$
The theorem can be stated by: $\majN\approx\majL$, where $\approx$ means that the two sets differ by a $\mu$-negligeable set.

The proof of this result uses stochastic methods which will be explained in the next section.


\section{Conditioning}

By Doob's h-processes, it is possible to condition the random walk to exit at a fixed point $\theta\in \partial S$ (see \cite{Doo57} and \cite{Dyn69}). The probability $\mathbb{P}_z^\theta$ on $\Omega$ thus obtained satisfies a strong Markow property and one has the following property:

\begin{proposition}    \label{prop:desintegration}
Let $F$ be a non-negative random variable on $\Omega$. Then
$$\mathbb{E}_z[F]=\intop_{\partial S} \mathbb{E}_z^\theta [F] d\mu_z (\theta).$$
\end{proposition}

The probability $\mathbb{P}_z^\theta$ satisfies an asymptotic zero-one law: if an event $A$ is asymptotic (\textit{i.e.} if it is invariant under the shift operator $\Theta$) then, for all $\theta\in\partial S$, the map $z\mapsto \mathbb{P}_z^\theta (A)$ is constant on $S$ and equals either $0$ or $1$.
The reader can refer to \cite{Mou94}, \cite{Bro78} and \cite{Dur84} for more details.

As mentioned above, we shall use probabilistic methods. We shall therefore define stochastic analogues of non-tangential convergence and boundedness notions.
Let $u$ be a harmonic function. Let $\widetilde{\majN^{**}}$ be the set of trajectories $\omega$ such that $|u|$ is bounded on the thickened trajectory
 $\{ y\in S| d(y,\omega)\leq m_1\}$:

$$\widetilde{\majN^{**}}=\{\omega\in\Omega \, | \, \widetilde{N^*}(\omega) <+\infty \} \text{ where } \widetilde{N^*}(\omega)=\sup \{ |u(y)| \, | \, y\in S, d(y,\omega)\leq m_1 \}.$$
We also define the set
$$\majL^{**}=\{\omega\in\Omega \, | \, \lim_{n\to\infty} u(X_n(\omega)) \text{ exists and is finite} \}.$$

These two events are asymptotic, so by asymptotic zero-one law, quantities $\mathbb{P}_z^\theta(\widetilde{\majN^{**}})$ and $\mathbb{P}_z^\theta(\majL^{**})$ have values $0$ or $1$ and do not depend on $z$. We thus define the sets

$$\widetilde{\majN^*}=\{\theta\in\partial S \, | \, \mathbb{P}_o^\theta (\widetilde{\majN^{**}})=1 \}
\text{ and } \majL^*=\{\theta\in\partial S \, | \, \mathbb{P}_o^\theta (\majL^{**})=1 \}.$$

We say that $u$ is \textit{stochastically bounded} at $\theta\in\partial S$ if $\theta\in \widetilde{\majN^*}$ and that $u$ \textit{converges stochastically} at $\theta$ if $\theta\in \majL^*$. For every $r\in \mathbb{R}$, the event $\{ \omega \, | \, \lim_{n\to\infty} u(X_n(\omega))\leq r \}$ is asymptotic, thus if $\theta\in\majL^*$, $\lim_{n\to\infty} u(X_n)$ is $\mathbb{P}_o^\theta$-almost surely constant.

We now prove a stochastic analogue of our main result.

\begin{proposition}    \label{prop:stochastic result}
Given a harmonic function $u$, one has the $\mu$-almost inclusion
$$\widetilde{\majN^{*}}\simsubset \majL^{*}.$$
\end{proposition}

\begin{proof}

We will first prove the $\mathbb{P}_o$-almost inclusion $\widetilde{\mathcal{N}^{**}}\simsubset\mathcal{L}^{**}$. For $m\in\mathbb{N}$, denote by $\widetilde{\mathcal{N}_m^{**}}$ the set of trajectories $\omega$ such that $|u|$ is bounded by $m$ on the thickened trajectory ${\{ y\in S \, | \, d(y,\omega)\leq m_1 \}}$. By countable union, it is sufficient to prove that for all $m$, ${\widetilde{\mathcal{N}_{m}^{**}}\simsubset\mathcal{L}^{**}}$.
Denote by $T_m$ the stopping time 
$$ T_m:=\inf \{ n\geq 0 \, | \, \max \{|u(y)| \, | \, y\in S, d(y,X_n)\leq m_1 \} >m \} .$$
Remark that $\widetilde{\mathcal{N}_m^{**}}=\{ T_m=+\infty\}$.
Since $u$ is harmonic, $(u(X_{n}))_{n\in\mathbb{N}}$
is a martingale for the probability $\mathbb{P}_{o}$ and thus $(u(X_{n\wedge T_m}))_{n}$
is a martingale too.
With our choice of stopping time $T_m$, for all $n\in\mathbb{N}$, $|u(X_{n\wedge T_m})| \leq \max \{ m, |u(X_o)| \}$,
which implies by the martingale theorem that the stopped martingale converges $\mathbb{P}_{o}$-almost surely.
In particular, $(u(X_{n}))_{n}$ converges $\mathbb{P}_o$-almost surely on the event
$\widetilde{\mathcal{N}_m^{**}}$.
We thus proved that $\widetilde{\mathcal{N}_{m}^{**}}\simsubset\mathcal{L}^{**}$
and since $\widetilde{\mathcal{N}^{**}}=\bigcup_m \widetilde{\mathcal{N}_{m}^{**}}$ we obtain that $\widetilde{\mathcal{N}^{**}}\simsubset\mathcal{L}^{**}$.

Using proposition \ref{prop:desintegration}, we have 
\[
0=\mathbb{P}_{o}(\widetilde{\mathcal{N}^{**}}\setminus\mathcal{L}^{**})=\intop_{\partial S} \mathbb{P}_{o}^{\theta}(\widetilde{\mathcal{N}^{**}}\setminus\mathcal{L}^{**})d\mu_{o}(\theta).
\]
 Then, $\mathbb{P}_{o}^{\theta}(\widetilde{\mathcal{N}^{**}}\setminus\mathcal{L}^{**})=0$
for $\mu$-almost all $\theta\in\partial S$ and $\widetilde{\mathcal{N}^{*}}\simsubset\mathcal{L}^{*}$. 
\end{proof}

We end this section with the case of bounded harmonic functions (\cite{Woe00}, \cite{Anc90}).

\begin{lemma}	\label{lem:Representation of bounded harmonic functions}
A bounded harmonic function $u$ on $S$ converges non-tangentially and stochastically for $\mu$-almost all point $\theta\in\partial S$ and the unique function $f\in L^\infty(\partial S,\mu)$ such that 
$$u(x)=\intop_{\partial S} f(\theta )d\mu_x (\theta )=\mathbb{E}_x[ f(X_{\infty})]$$
is $\mu$-a.e. the non-tangential and stochastic limit of $u$.
\end{lemma}


\section{Geometric lemmas}

We begin this section by showing that a hyperbolicity inequality holds for points on the boundary $\partial S$: for all $x,y,z\in S\cup\partial S$,
\begin{equation}
 (x,y)_o\geq \min \{ (x,z)_o,(y,z)_o \} -2\delta.  \label{gromov hyberbolicity bis}
\end{equation}

To see this, choose, for $\epsilon>0$, sequences in $S$ with $x_i\to x$, $y_i\to y$, $z_i\to z$ and $z_i'\to z$ such that $\liminf_{i,j} (x_i,z_j)_o \geq (x,z)_o-\epsilon$ and $\liminf_{i,j} (z_i',y_j)_o \geq (z,y)_o-\epsilon$. Then, take $\liminf_{i,j}$ through $(x_i,y_j)_o\geq \min \{ (x_i,z_j)_o , (z_j,z_i')_o , (z_i',y_j)_o \} -2\delta$ (note that $\liminf_{i,j} (z_j,z_i')_o=+\infty$).

We will also need the fact that for all $x\in S$, $\xi\in\partial S$, and all geodesic ray $\gamma$ from $o$ to $\xi$,
\begin{equation}
 d(x,\gamma)-2\delta \leq (o,\xi)_x \leq d(x,\gamma)+2\delta. \label{distance geodesic ray}
\end{equation}

By inequality (\ref{thin triangles}), for all $i$, $d(x,\gamma($\textlbrackdbl $0,i$ \textrbrackdbl $ ))-2\delta \leq (o,\gamma(i))_x \leq d(x,\gamma($\textlbrackdbl $0,i$ \textrbrackdbl $ ))$.
Since $d(x,\gamma(i))\to\infty$, for $i$ large enough, $d(x,\gamma($\textlbrackdbl $0,i$ \textrbrackdbl $ ))=d(x,\gamma)$ and hence for $i$ large enough,
$$d(x,\gamma)-2\delta\leq (o,\gamma(i))_x \leq d(x,\gamma).$$
Combining this inequality with the fact that if $(\xi_i)_i$ is a sequence such that $\xi_i \to \xi$, thus $(o,\xi)_x -2\delta \leq \liminf_i (o,\xi_i)_x \leq (o,\xi)_x$ (see \cite{Bri99}), we obtain inequality (\ref{distance geodesic ray}).

Using the hyperbolicity of the graph, we shall prove two lemmas, and deduce three corollaries. In order to prove one of these lemmas (lemma \ref{lem:geometric lemma1}), we shall need some Harnack inequalities (see \cite{Anc88} and \cite{Anc90}):

\begin{theorem}[Harnack inequality]
\label{thm:Harnack on balls}
 Let $u$ be a non-negative superharmonic function. For all $x,y\in S$,
$$\left(\frac{c_0}{\ell}\right)^{d(x,y)} u(y) \leq u(x) \leq \left(\frac{\ell}{c_0}\right)^{d(x,y)}u(y).$$
\end{theorem}

The following theorem is a version of the so-called Harnack inequality at infinity of A.~Ancona.

\begin{theorem}[submultiplicativity of the Green function]
\label{thm:submultiplicativity Green function}
For any $r>0$, there exists a constant $C=C(r)$ such that, for all $x,z\in S$ and all $y\in S$ at distance at most $r$ from any geodesic segment between $x$ and $z$,
$$G(x,z)\leq C \cdot G(x,y)G(y,z).$$
\end{theorem}

We can now state the geometric lemmas. They will be of central importance in the proof of the main result (theorem \ref{thm:main result}).

\begin{lemma}    \label{lem:geometric lemma1}
Given $\alpha>0$, there exists a constant $C>0$ such that,
for all point $x\in S$ and all $\theta \in \partial S$, \[
\mu_{x}(\{\xi\in\partial S \, | \, (\xi,\theta)_{x}\geq\alpha\})\geq C.\]
\end{lemma}

\begin{proof}

First, we will show that there exists $N=N(\alpha)>0$ such that for all $x\in S$, all $\theta\in\partial S$ and all $y$ on a geodesic from $x$ to $\theta$ with $d(x,y)\geq N$, we have
$$\mu_y( A_{x,\alpha}^\theta ) > \frac{1}{2}$$
where
$$A_{x,\alpha}^\theta := \{ \xi\in\partial S \, | \, (\xi,\theta)_x \geq \alpha \}.$$
We have 
$$\mu_y(\partial S \setminus A_{x,\alpha}^\theta)=\intop_{\partial S \setminus A_{x,\alpha}^\theta} \frac{d\mu_y}{d\mu_x}(\xi) d\mu_x(\xi)$$
and
$$\frac{d\mu_y}{d\mu_x}(\xi)=\lim_{z\to \xi} \frac{G(y,z)}{G(x,z)}.$$
Let $\xi\in \partial S \setminus A_{x,\alpha}^\theta$. Denote by $\gamma$ a geodesic ray from $x$ to $\xi$, by $x_j=\gamma(4j\delta)$ and by $V_j$ the closure in $S\cup\partial S$ of $\{ z\in S \, | \, (z,x_j)_x > d(x,x_j)-3\delta \}$. 

We claim that there exist $j,N_1\in\mathbb{N}$ depending only on $\alpha$ such that if $y$ is on a geodesic from $x$ to $\theta$ and $d(x,y)\geq N_1$, then $y\notin V_j$ (see figure \ref{figure lem:geometric lemma1}). Indeed, choose $j$ such that $d(x,x_j)-3\delta=4j\delta -3\delta >\alpha+4\delta$. By the hyperbolicity inequality (\ref{gromov hyberbolicity bis}),

\begin{figure}
\begin{center}

\includegraphics[width=4.5cm,height=4.5cm]{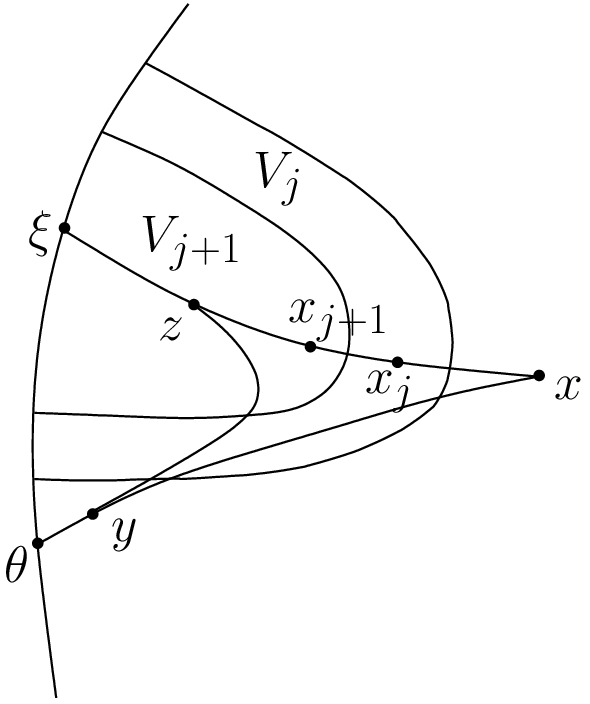}

\caption{Lemma \ref{lem:geometric lemma1}}
\label{figure lem:geometric lemma1}
\end{center}
\end{figure}

$$\alpha > (\xi,\theta)_x\geq \min \{ (\xi,y)_x,(y,\theta)_x \} -2\delta$$
and if $y$ is on a geodesic ray from $x$ to $\theta$, there exists $N_1$ depending only on $\alpha$ such that $d(x,y)\geq N_1$ implies $(y,\theta)_x >\alpha +2\delta$. Thus for such a $y$, $(\xi,y)_x\leq \alpha+2\delta$. Using once again the hyperbolicity inequality,
$$\alpha+2\delta \geq (\xi,y)_x\geq \min \{ (\xi,x_j)_x,(x_j,y)_x \}-2\delta.$$
Since $\xi\in V_j$, $(\xi,x_j)_x \geq d(x,x_j)-3\delta > \alpha+4\delta$, thus $(x_j,y)_x\leq \alpha+4\delta$ and $y\notin V_j$, which proves the claim.

Now, we can verify that for all $z\in V_{j+1}$ and all $y\notin V_j$, the distance between $x_j$ and a geodesic segment between $y$ and $z$ is at most $50\delta$ (\cite{Anc90} p85). Then by theorem \ref{thm:submultiplicativity Green function}, $G(y,z)\leq C_1(\delta) \cdot G(y,x_j)G(x_j,z)$. We can now apply theorem \ref{thm:Harnack on balls} to $G(\cdot,z)$ and $G(y,\cdot)$ ($G(y,\cdot)$ is superharmonic for the admissible function $p^*$) and we obtain $G(y,z)\leq C(\delta, \alpha) \cdot G(y,x)G(x,z)$. Making $z\to \xi$, $z\in V_{j+1}$, we obtain
$$\frac{d\mu_y}{d\mu_x}(\xi) \leq C\cdot G(y,x).$$
Thus, $\mu_y(\partial S \setminus A_{x,\alpha}^\theta)\leq C \cdot G(y,x)$. Since $G$ has a uniform exponential decay at infinity (proposition \ref{prop:exponential decay}), there exists $N$ depending only on $\alpha$ such that for all $y$ on a geodesic from $x$ to $\theta$ with $d(x,y)\geq N$,
$$\mu_y( A_{x,\alpha}^\theta ) > \frac{1}{2}.$$
By Harnack inequality (theorem \ref{thm:Harnack on balls}), if $x,y\in S$ with $d(x,y)=N$, 
$$\mu_x( A_{x,\alpha}^\theta ) \geq \left(\frac{c_0}{\ell} \right)^N \mu_y( A_{x,\alpha}^\theta )\geq C>0.$$

\end{proof}

For a borelian set $E\subset\partial S$, we denote $\Gamma_c(E) :=\bigcup_{\theta\in E} \Gamma_c^\theta$.

\begin{lemma}   \label{lem:geometric lemma2}
There exists $\eta>0$ and $c_1>0$ such that for all $c > c_1$ and all borelian
sets $E\subset\partial S,$ one has\[
\forall x\not\in\Gamma_{c}(E),\, \mathbb{P}_{x}(X_{\infty}\not\in E)\geq\eta.\]
\end{lemma}

\begin{proof}
Recall that $K$ is the constant provided by the assumption (GRBD). Let $c_1=K+6\delta$.
Fix $c>c_1$, a borelian set $E$ in $\partial S$ and $x\not\in\Gamma_{c}(E)$. 
Choose a geodesic ray $\bar\gamma$ from $o$ to a point $\bar\xi\in\partial S$ such that $d(x,\bar\gamma)\leq K$.

\begin{figure}
\begin{center}

\includegraphics[width=4.5cm,height=5cm]{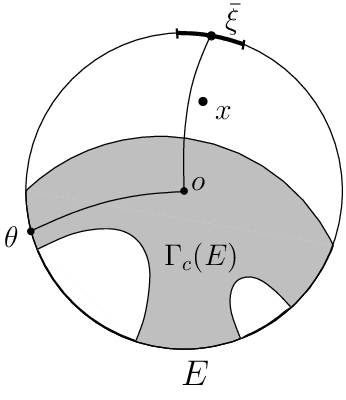}

\caption{Lemma \ref{lem:geometric lemma2}}
\label{figure lem:geometric lemma2}
\end{center}
\end{figure}

In order to use lemma \ref{lem:geometric lemma1}, we will show that there exists a constant $\alpha>0$ depending only on $\delta$ and $K$ such that $\{\xi\in\partial S \, | \, (\xi,\bar\xi)_{x}\geq\alpha\} \subset \partial S\setminus E$.

For $\theta\in E$, we want to bound uniformly from above the product $(\theta,\bar\xi)_x$. Inequality (\ref{gromov hyberbolicity bis}) gives
$$\min \{ (\theta,\bar\xi)_x , (o,\theta)_x \} \leq (o,\bar\xi)_x +2\delta.$$

By inequality (\ref{distance geodesic ray}), $( o ,\bar\xi )_x \leq d(x,\bar\gamma)+2\delta \leq K+2\delta$, so $\min \{ (\theta,\bar\xi)_x , (o,\theta)_x \} \leq K+4\delta$. Again by inequality (\ref{distance geodesic ray}), denoting by $\gamma$ a geodesic ray from $o$ to $\theta$, $(o,\theta)_x \geq d(x,\gamma)-2\delta \geq c-2\delta >K+4\delta$ and  $\min \{ (\theta,\bar\xi)_x , (o,\theta)_x \} = (\theta,\bar\xi)_x \leq K+4\delta$.

Therefore, $\{\xi\in\partial S \, | \, (\xi,\bar\xi)_x\geq K+5\delta \}\cap E=\emptyset$. By lemma \ref{lem:geometric lemma1}, there exists $\eta>0$ depending only on $\delta$ such that 
$$\mathbb{P}_x(X_\infty \not\in E)\geq \eta.$$

\end{proof}

\begin{corollary}   \label{cor:the end lemma}
Let $E$ be a borelian set of $\partial S$, $x\in S$ and $c>c_1$.
For $\mu$-almost all $\theta\in E$, $\mathbb{P}_{x}^{\theta}$-a.s., the random
walk ``ends in $\Gamma_{c}(E)$'' (Formally, for $\mathbb{P}_{x}^{\theta}$-almost all $\omega$, there exists $N\in\mathbb{N}$ such that for all $n\geq N$, $X_{n}(\omega ) \in\Gamma_{c}(E)$). 
\end{corollary}

\begin{proof}
Let $f_{E}(x):=\mathbb{P}_{x}(X_{\infty}\in E)=\mathbb{E}_{x}[\ind{E}(X_{\infty})]$.
As a consequence of the representation lemma \ref{lem:Representation of bounded harmonic functions} of bounded harmonic
functions, for $\mu$-almost all $\theta$, one has 
\[
\forall x\in S, \, \mathbb{P}_{x}^{\theta}[\lim_{n\to\infty}f_{E}(X_{n})=\ind{E}(\theta)]=1.\]
Because of lemma \ref{lem:geometric lemma2}, there exists $\eta>0$ such that \[
\forall x\not\in\Gamma_{c}(E),\, f_{E}(x)\leq1-\eta.\]
Thus for all $x\in S$ and for $\mu$-almost all $\theta\in E$, $\mathbb{P}_{x}^{\theta}$-a.s.,
$X_{n}$ is in $\Gamma_{c}(E)$ for $n$ large enough.
\end{proof}

Given a tube $\Gamma_c^\theta$ and $R>0$, the set $\Gamma_c^\theta \setminus B(o,R)$ is called a \textit{spike} of $\Gamma_c^\theta$.

\begin{corollary}   \label{cor:spikes lemma}
Let $c>c_1$ and $E$ be a borelian set of $\partial S$. Then, for all
$\theta\in\partial S$ such that ${\displaystyle \limNT{x}{\theta} \mathbb{P}_{x}(X_{\infty}\in E)=1}$,
$\Gamma_{c}(E)$ contains spikes of every tube with $\theta$ as vertex.

In particular, it is the case for $\mu$-almost all $\theta\in E$ by the
bounded harmonic function representation lemma \ref{lem:Representation of bounded harmonic functions}.
\end{corollary}

\begin{proof}
Fix $\theta\in\partial S$ such that ${\displaystyle \limNT{x}{\theta} \mathbb{P}_{x}(X_{\infty}\in E) = \limNT{x}{\theta}f_{E}(x)=1}$ and let $\Gamma_{e}^{\theta}$ be a tube with vertex $\theta$. 
By contradiction,
assume that $\Gamma_{c}(E)$ does not contain any spike of this tube.
Then, for each $R>0$, there exists $x\in\Gamma_{e}^{\ \theta}\setminus\Gamma_{c}(E)$
such that $d(o,x)>R$.

Let $(x_{k})_{k\in\mathbb{N}}$ be a sequence in $\Gamma_{e}^{\ \theta}\setminus\Gamma_{c}(E)$
such that $d(o,x_{k})>k$. Then $(x_{k})_{k}$
converges to $\theta$ staying in $\Gamma_e^\theta$ and we have $\lim_{k\to\infty} f_{E}(x_{k})=1.$
Since $x_{k}\not\in\Gamma_{c}(E)$, by lemma \ref{lem:geometric lemma2}, $f_{E}(x_{k})\leq1-\eta$,
a contradiction.
\end{proof}

The following corollary will not intervene later. However, it is remarkable to observe that the behaviour of a harmonic function on a tube $\Gamma_{c_0}^\theta$ controls the behaviour of this function on every tube $\Gamma_c^\theta$, $c>0$.

\begin{corollary}  
Given a harmonic function $u$, for all real $c>c_1$ one has $\mathcal{N}_{c}\approx\mathcal{N}$.
\end{corollary}

\begin{proof}
By definition, $\mathcal{N}\subset\mathcal{N}_{c}$. It is thus sufficient to show that $\mathcal{N}_{c}\simsubset\mathcal{N}$ for
$c>c_1$.
Let $c>c_1$ and denote by
$$A_c^m=\{\theta\in\partial S \, | \, N_{c}^{\theta}(u)\leq m\}$$
the set of points $\theta\in\partial S$ such that $u$ is bounded by $m$ on $\Gamma_c^\theta$.
As $\mathcal{N}_{c}$ is the countable union of the $A_c^m$, we need only to show that $A_c^m\simsubset\mathcal{N}$ for all $m$.
By definition of $A_c^m$, $|u|$ is bounded by $m$ on $\Gamma_{c}(A_c^m)$.
Using corollary \ref{cor:spikes lemma}, we obtain that for $\mu$-almost all
points $\theta\in A_c^m$, $\Gamma_{c}(A_c^m)$ contains spikes
of all tube with $\theta$ as vertex. On these spikes, the function $u$ is bounded,
and therefore, by local finiteness, $u$ is bounded on the tubes,
which means that $\theta$ is in $\mathcal{N}$. Finally, $A_c^m\simsubset\mathcal{N}$ and therefore, 
$$\mathcal{N}_{c}\simsubset\mathcal{N}.$$
\end{proof}


\section{proof of the main result}

With geometric lemmas and the stochastic result (proposition \ref{prop:stochastic result}) in hand, we can now prove theorem \ref{thm:main result}.

\begin{proof}
As above, let us denote by
$$A_c^m=\{\theta\in\partial S \, | \, N_c^\theta (u)\leq m\}.$$

Since $\majN_c$ is a countable union of the sets $A_c^m$, it is sufficient to prove that for all $m$ and all $c>c_1 +m_1$, $A_c^m\simsubset \majL_{c-m_1}$. Then, we will have $\majN_c\simsubset \majL_{c-m_1}$ and since for $c>c'$, $\majL_{c}\subset \majL_{c'}$, we can conclude that $\majN\simsubset \majL$.

Let $c>c_1+m_1$.
We shall first prove that $A_c^m\simsubset \majL^*$. Applying corollary \ref{cor:the end lemma} to the borelian set $A_c^m$, we get: for $\mu$-almost all point $\theta\in A_c^m$, $\mathbb{P}_z^\theta$-almost surely, $(X_k)_{k\geq 0}$  ends in $\Gamma:=\Gamma_{c-m_1} (A_c^m)$. Let $\theta$ be such a point. The key point is that for all $x\in \Gamma$ and all $y\in S$ such that $d(x,y)\leq m_1$, $|u(y)|\leq m$. It implies in particular that for $\mathbb{P}_z^\theta$-almost all $\omega$, there exists $N\in\mathbb{N}$ such that for all $n\geq N$ and all $y\in S$ such that $d(y,X_n(\omega))\leq m_1$, we have $|u(y)|\leq m$. By local finiteness, $\mathbb{P}_z^\theta$-almost surely,  $\widetilde{N^*}<+\infty$. Thus, $\mathbb{P}_z^\theta(\widetilde{\majN^{**}})=1$, $\theta \in \widetilde{\majN^*}$ and hence $A_c^m \simsubset\widetilde{\majN^*}$. However, by proposition \ref{prop:stochastic result}, $\widetilde{\majN^*}\simsubset \majL^*$, so 
$$A_c^m\simsubset \majL^*.$$

Let us now prove that $A_c^m\simsubset \majL_{c-m_1}$.  As shown above, for $\mu$-almost all $\theta\in A_c^m$, $\mathbb{P}_z^\theta$-almost surely, $(u(X_n))_n$ has a finite limit $\ell(\theta )$. It defines a function $\ell$ on $A_c^m$. We use again corollary \ref{cor:the end lemma}: for $\mu$-almost all $\theta\in A_c^m$, $\mathbb{P}_z^\theta$-almost surely, $X_n$ is in $\Gamma$ for $n$ big enough. This together with the fact that $|u|$ is bounded by $m$ on $\Gamma$ implies that $| \ell | \leq m$ on $A_c^m$.

We will conclude this proof using a method of J.~Brossard \cite{Bro78}. We will decompose $u$ on $\Gamma$ as a sum of three functions which will have non-tangential limits at almost all points of $A_c^m$.

We define the function 
$$f(z):=\mathbb{E}_z \left[(\ell \cdot \ind{ A_c^m } )(X_{\infty})\right].$$
By the representation lemma \ref{lem:Representation of bounded harmonic functions}, $f$ is a bounded harmonic function which converges non-tangentially at $\mu$-almost all point $\theta\in\partial S$ to $(\ell \cdot \ind{A_c^m})(\theta )$. Denote by $\tau$ the exit time of the set $\Gamma$ and $\tau_k$ the exit time of $B(o,k)$. Since $u$ is bounded and harmonic on the thickened set $\widetilde{\Gamma}\cap \widetilde{B(o,k)}=\{y\in S \, | \, d(y,\Gamma\cap B(o,k))\leq m_1\}$, which is a bounded set,
${u(z)=\mathbb{E}_z[u(X_{\tau\wedge \tau_k})]}$. If $\tau=+\infty$, $\mathbb{P}_z$-almost surely, $(X_n)_n$ converges to a point $X_\infty \in A_c^m$, so $\mathbb{P}_z$-almost surely, $(u(X_n))_n$ goes to $\ell(X_\infty )$. We can extend $u$ to $A_c^m$ by setting $u(\theta ):=\ell(\theta )$ for $\theta\in A_c^m$. Since $|u|$ is bounded by $m$ on $\Gamma$, we can apply Lebesgue's theorem to obtain
$$\forall z\in\Gamma, u(z)=\mathbb{E}_z[u(X_\tau )].$$
Decomposing the event $\{ X_\infty\in A_c^m \}$ into the union ${\{\tau <\infty ; X_\infty\in A_c^m \}\cup \{\tau =\infty \}}$ we obtain, for $z\in\Gamma$,
\begin{align*}
u(z) &= \mathbb{E}_z[u(X_\tau ) \cdot \ind{ \{\tau <\infty\} }]+\mathbb{E}_z[u(X_\infty ) \cdot \ind{\{\tau =\infty\} }] \\
     &= \mathbb{E}_z[u(X_\tau ) \cdot \ind{\{\tau <\infty\} }]+\mathbb{E}_z[u(X_\infty ) \cdot \ind{\{ X_\infty \in A_c^m\} }]\\
     &  -\mathbb{E}_z[u(X_\infty ) \cdot \ind{\{ X_\infty \in A_c^m\} } \cdot \ind{\{ \tau <\infty\} }].
\end{align*}
This is exactly the announced decomposition. Indeed, denoting $g(z)=\mathbb{E}_z[u(X_\tau ) \cdot \ind{\{\tau <\infty\} }]$ and $h(z)=-\mathbb{E}_z[u(X_\infty ) \cdot \ind{\{ X_\infty \in A_c^m\} }.\ind{\{\tau <\infty\} }]$, we have $u=f+g+h$ on $\Gamma$.

It remains to prove that at almost all point $\theta$ in $A_c^m$, the functions $g(z)$ and $h(z)$ converge to zero when $z$ goes to $\theta$ staying in the tube $\Gamma_{c-m_1}^\theta$.
Since $u$ is bounded on $\widetilde{\Gamma}=\{y\in S \, | \, d(y,\Gamma)\leq m_1\}$, if $\tau<\infty$, then $|u(X_\tau )|\leq m$ and obviously 
$$|g(z)|\leq m \cdot \mathbb{P}_z(\tau<\infty).$$
In the same way, for almost all $\theta\in A_c^m$, $|u(\theta )|\leq m$, so we obtain easily by conditioning
$$|h(z)|\leq m \cdot \mathbb{P}_z(\tau<\infty).$$
It is now sufficient to prove that for almost all $\theta \in A_c^m$, $\mathbb{P}_z(\tau<\infty)$ goes to zero when $z$ goes to $\theta$ staying in $\Gamma_{c-m_1}^\theta$. This follows from lemma \ref{lem:geometric lemma2}.
Indeed, there exists $\eta>0$ such that 
$$\forall z\not\in \Gamma, \mathbb{P}_z(X_\infty \not\in A_c^m )\geq \eta$$
and in particular this holds for all $z\in \widetilde\Gamma\setminus\Gamma$.
The strong Markov property implies that for all $z\in\Gamma$,
\begin{align*}
\mathbb{P}_z(X_\infty \not\in A_c^m ) &= \mathbb{P}_z(\{ X_\infty \not\in A_c^m \}\cap \{ \tau<\infty \} ) \\
	 &= \sum_{i=1}^\infty \mathbb{P}_z(\{X_\infty \not\in A_c^m\}\cap \{ \tau=i \} ) \\
	 &= \sum_{i=1}^\infty \mathbb{E}_z\left[ \mathbb{E}_z \left[ \ind{\{X_\infty \not\in A_c^m\}} \cdot \ind{ \{ \tau=i \}} | \mathcal{F}_i \right] \right] \\
	 &= \sum_{i=1}^\infty \mathbb{E}_z\left[ \mathbb{E}_z \left[ \ind{\{X_\infty \not\in A_c^m\}}  | \mathcal{F}_i \right]\cdot \ind{ \{ \tau=i \}}  \right] \\
	&= \sum_{i=1}^\infty \mathbb{E}_z[ \mathbb{P}_{X_i}(X_\infty \not\in A_c^m) \cdot \ind{\{\tau=i\} }] \\
	&= \mathbb{E}_z[ \mathbb{P}_{X_\tau}(X_\infty \not\in A_c^m ) \cdot \ind{\{ \tau<\infty\} }] \\
	&\geq \eta \cdot \mathbb{P}_z(\tau<\infty ).
\end{align*}
By lemma \ref{lem:Representation of bounded harmonic functions}, for almost all $\theta\in A_c^m$, $\limNT{z}{\theta} \mathbb{P}_z(X_\infty \not\in A_c^m )=0$. It follows that for almost all $\theta \in A_c^m$, $\mathbb{P}_z(\tau<\infty)$ goes to zero when $z$ goes to $\theta$ staying in $\Gamma_{c-m_1}^\theta$.\\
We thus have that $A_c^m \simsubset \majL_{c-m_1}$ and the theorem is proved.										 

\end{proof}


\end{document}